\documentclass[11pt]{article}
\usepackage{amsfonts,amsmath,amssymb,hyperref,accents,amsthm}
\usepackage{tikz}
\usetikzlibrary{calc,arrows,positioning,fit,petri}

\newcommand{\file}{}

\newcommand{\dis}{\displaystyle}
\newcommand{\txt}{\textstyle}

\newcommand{\noi}{\noindent}
\newcommand{\halmos}{\rule{1ex}{1.4ex}}
\newcommand{\QED}{\nopagebreak{\hspace*{\fill}$\halmos$\medskip}}

\newcommand{\quand}{\quad\mbox{and}\quad}

\newtheoremstyle{mythm}
  {}
  {}
  {\itshape}
  {}
  {\bfseries}
  {}
  {.5em}
  {#1 #2 \thmnote{(#3)}}

\theoremstyle{mythm}
\newtheorem{theorem}{Theorem}
\newtheorem{proposition}[theorem]{Proposition}
\newtheorem{lemma}[theorem]{Lemma}
\newtheorem{exercise}[theorem]{Exercise}
\newtheorem{corollary}[theorem]{Corollary}
\newtheorem{conjecture}[theorem]{Conjecture}

\newtheorem{counterex}[theorem]{Counterexample}

\newcommand{\bt}{\begin{theorem}}
\newcommand{\et}{\end{theorem}}
\newcommand{\bl}{\begin{lemma}}
\newcommand{\el}{\end{lemma}}
\newcommand{\bp}{\begin{proposition}}
\newcommand{\ep}{\end{proposition}}
\newcommand{\bcor}{\begin{corollary}}
\newcommand{\ecor}{\end{corollary}}
\newcommand{\br}{\begin{remark}\rm}
\newcommand{\er}{\end{remark}}
\newcommand{\bcon}{\begin{conjecture}}
\newcommand{\econ}{\end{conjecture}}
\newcommand{\bex}{\begin{exercise}}
\newcommand{\eex}{\end{exercise}}
\newcommand{\bcou}{\begin{counterex}}
\newcommand{\ecou}{\end{counterex}}

\newenvironment{Proof}[1][]{\noi\textbf{Proof #1}}{\QED}
\newcommand{\bpro}{\begin{Proof}}
\newcommand{\epro}{\end{Proof}}

\newcommand{\be}{\begin{equation}}
\newcommand{\ee}{\end{equation}}
\newcommand{\ba}{\begin{array}}
\newcommand{\ea}{\end{array}}
\newcommand{\bc}{\be\begin{array}{r@{\,}c@{\,}l}}
\newcommand{\ec}{\end{array}\ee}
\newcommand{\bac}{\begin{array}{r@{\,}c@{\,}l}}

\newcommand{\la}{\lambda}

\newcommand{\sig}{\sigma}

\newcommand{\si}{\ensuremath{\sigma}}

\newcommand{\Pc}{{\cal P}}

\newcommand{\R}{{\mathbb R}}
\newcommand{\N}{{\mathbb N}}
\newcommand{\Z}{{\mathbb Z}}

\renewcommand{\P}{{\mathbb P}}
\newcommand{\E}{{\mathbb E}}

\newcommand{\up}{\uparrow}
\newcommand{\down}{\downarrow}
\newcommand{\sub}{\subset}
\newcommand{\beh}{\backslash}

\newcommand{\asto}[1]{\underset{{#1}\to\infty}{\longrightarrow}}

\newcommand{\ti}{\tilde}

\newcommand{\ov}{\overline}

\newcommand{\di}{\mathrm{d}}
\newcommand{\half}{{[0,\infty)}}

\newcommand{\Lin}{L}
\newcommand{\Fo}{\ov F}
\newcommand{\No}{\ov N}
\newcommand{\Yo}{\ov Y}

\begin{document}


\renewcommand{\labelenumi}{{\rm(\roman{enumi})}}

\title{\vspace{-3cm}A simple rank-based Markov chain with self-organized
  criticality}
\author{Jan M. Swart\vspace{6pt}\\
{\small \' UTIA}\\
{\small Pod vod\'arenskou v\v e\v z\' i 4}\\
{\small 18208 Praha 8}\\
{\small Czech Republic}\\
{\small e-mail: swart@utia.cas.cz}
\vspace{4pt}}
\date{{\scriptsize\file}\quad\today}
\maketitle\vspace{-.7cm}

\begin{abstract}\noi
We introduce a self-reinforced point processes on the unit interval that
appears to exhibit self-organized criticality, somewhat reminiscent of the
well-known Bak-Sneppen model. The process takes values in the finite
subsets of the unit interval and evolves according to the following rules. In
each time step, a particle is added at a uniformly chosen position,
independent of the particles that are already present. If there are any
particles to the left of the newly arrived particle, then the left-most of
these is removed. We show that all particles arriving to the left of $p_{\rm
  c}=1-e^{-1}$ are a.s.\ eventually removed, while for large enough time,
particles arriving to the right of $p_{\rm c}$ stay in the system forever.
\end{abstract}

\vspace{.4cm}
\noi
%
{\it Acknowledgments.} Work sponsored by GA\v{C}R grants: P201/12/2613 and
15-08819S.


{\setlength{\parskip}{-2pt}\tableofcontents}
\newpage

\section{Introduction and results}

\subsection{Main results}

Let $(U_k)_{k\geq 1}$ be an i.i.d.\ collection of uniformly distributed
$[0,1]$-valued random variables. For each finite subset $x$ of $[0,1]$, we
inductively define a sequence $X^x=(X^x_k)_{k\geq 0}$ of random finite subsets
of $[0,1]$ by $X^x_0:=x$, $M^x_{k-1}:=\min(X^x_{k-1}\cup\{1\})$ and
\be\label{induct}
X^x_k:=\left\{\ba{ll}
X^x_{k-1}\cup\{U_k\}\quad&\mbox{if }U_k<M^x_{k-1},\\[5pt]
(X^x_{k-1}\cup\{U_k\})\beh\{M^x_{k-1}\}\quad&\mbox{if }U_k>M^x_{k-1}.
\ea\right.\quad(k\geq 1).
\ee
In words, this says that $M^x_{k-1}$ is the minimal element of $X^x_{k-1}$ and
that the set $X^x_k$ is constructed from $X^x_{k-1}$ by adding $U_k$, and in
case that $M^x_{k-1}<U_k$, removing $M^x_{k-1}$ from $X^x_{k-1}$. Since the
$(U_k)_{k\geq 1}$ are i.i.d.\ and $X^x_k$ is a function of $X^x_{k-1}$ and
$U_k$, it is clear that $X^x$ is a Markov chain. (In fact, we have just given
a \emph{random mapping representation} for it.) The state space of $X^x$ is
the set $\Pc_{\rm fin}[0,1]$ of all finite subsets of $[0,1]$, which is
naturally isomorphic to the space of all simple counting measures on $[0,1]$
(i.e., $\N$-valued measures $\nu$ such that $\nu(\{u\})\leq 1$ for all
$u\in[0,1]$). We equip this space with the topology of weak convergence and
the associated Borel-\si-algebra.

The process $X^x$ is an example of a Markov process with self-reinforcement
(compare \cite{Pem07}), since the number of particles in the system can grow
without bound and influences the fate of newly arrived particles. As we will
see in a moment, it also appears to exhibit self-organized criticality in a
way that is reminiscent of the well-known Bak-Sneppen model.
The empirical distribution function $F^x_k(q):=\big|X^x_k\cap[0,q]\big|$ can
loosely be interpreted as the profile of a canyon being cut out by a river. If
$U_k<M^x_{k-1}$, then the river cuts deeper into the rock. If $U_k>M^x_{k-1}$,
then the slope of the canyon between $U_k$ and the river is eroded one step
down.

Our first result says that particles arriving on the left of the critical
point $p_{\rm c}:=1-e^{-1}$ are eventually removed from the system, but for
large enough time, particles arriving on the right of $p_{\rm c}$ stay in the
system forever.

\bt[A.s.\ behavior of the minimum]\label{T:limsup}
For any finite $x\sub[0,1]$, one has
\be
\limsup_{k\to\infty}M^x_k=1-e^{-1}\quad{\rm a.s.}
\ee
\et

To understand Theorem~\ref{T:limsup} better, note that for each $0\leq q\leq
1$, the restriction $X^x_k\cap[0,q]$ of $X^x_k$ to $[0,q]$ is a Markov
chain. Indeed, particles arriving on the right of $q$ just have the effect
that in each time step, with probability $1-q$, the minimal element of
$X^x_k\cap[0,q]$, if there is one, is removed, while no new particles are
added inside $[0,q]$. Theorem~\ref{T:limsup} says that this Markov chain is
recurrent for $q<p_{\rm c}$ and transient for $q>p_{\rm c}$. For any
$q\in[0,1]$, let
\be\label{tauemp}
\tau^q_\emptyset:=\inf\{k>0:X^\emptyset_k\cap[0,q]=\emptyset\}
\ee
be the first time the restricted process $X^\emptyset_k\cap[0,q]$ returns to
the empty set. Letting $\P^x$ denote the law of $X^x$, we have the following
result.

\bt[Recurrence versus transience]\label{T:reccur}
Let $p_{\rm c}:=1-e^{-1}$. Then
\be\ba{ll}\label{retn}
\dis\E^\emptyset[\tau^q_\emptyset]=\big(1+\log(1-q)\big)^{-1}
\qquad&\dis(q<p_{\rm c}),\\[5pt]
\dis\E^\emptyset[\tau^q_\emptyset]=\infty\quand
\P^\emptyset[\tau^q_\emptyset<\infty]=1\qquad&\dis(q=p_{\rm c}),\\[5pt]
\dis\P^\emptyset[\tau^q_\emptyset=\infty]>0\qquad&\dis(q>p_{\rm c}).
\ec
\et

Numerical simulations strongly suggest that at $q=p_{\rm c}$, the probability
$\P[\tau^q_\emptyset\geq k]$ decays as $k^{-1/2}$. We briefly comment on this
in Section~\ref{S:conclude}. In particular, it seems this can be proved
provided it can be shown that a certain martingale occurring in our proofs
behaves on long time scales like Brownian motion. Such a proof would establish
self-organized criticality for our process. Our process is self-organized in
the sense that it finds the transition point $p_{\rm c}$ by itself. In
particular, one does not have to tune a parameter of the model to exactly the
right value to see the (presumed) power-law critical behavior at~$p_{\rm c}$.

In the positive recurrent regime $q<p_{\rm c}$, it is not hard to show that
the process is ergodic, and as a result we also obtain the following
result. Below, we call a subset of $[0,p_{\rm c})$ \emph{locally finite} if
its intersection with any compact subset of $[0,p_{\rm c})$ is finite.

\bt[Ergodicity of restricted process]\label{T:ergod}
There exists a random, locally finite subset
$X_\infty\sub[0,p_{\rm c})$ such that, regardless of the initial state $x$,
\be\label{ergod}
\P\big[X^x_k\cap[0,q]\in\cdot\,\big]
\asto{k}\P\big[X_\infty\cap[0,q]\in\cdot\,\big]
\qquad(0<q<p_{\rm c}),
\ee
where $\to$ denotes convergence of probability measures in total variation
norm distance. The random point set $X_\infty$ a.s.\ consists of infinitely
many points.
\et

\subsection{Relation to known models}\label{S:other}

\subsubsection*{Models for email communication}

Our model appears to be new, but it is similar to a number of other models
that have been studied in the literature. To start with the simplest one,
consider the following model for email communication. (The model described
here is from \cite{FS16}, but very similar to a model introduced in
\cite{GC09}; see also \cite{Tok15} for a similar model with an interpretation
coming from mathematical finance.) A person receives emails at times of a
Poisson process with intensity $\la_{\rm in}$, and assigns to each email a
priority. Priorities are i.i.d.\ and uniformly distributed on $[-\la_{\rm
    in},0]$ (with $0$ the highest priority). At times of a Poisson process
with intensity $\la_{\rm out}:=1$, the person answers the email in the inbox
with the highest priority, if there is one, and does nothing
otherwise. Letting $X_t\sub[-\la_{\rm in},0]$ denote the set of priorities of
unanswered emails that are at time $t$ in the inbox, we observe that the
process is consistent in the sense that for each $0\leq\la\leq\la_{\rm in}$,
the process $(X_t\cap[-\la,0])_{t\geq 0}$ is a Markov process, that has the
same dynamics as the original model but with $\la_{\rm in}$ replaced by
$\la$. In this sense, the model is similar to our canyon model, but the email
model has the stronger property that if we throw away even more information
and only consider the number $N^\la_t:=|X_t\cap[-\la,0]|$ of emails with
priority above $-\la$, then even this is a Markov process. Indeed, it is easy
to see that $N^\la_t$ jumps $n\mapsto n+1$ with rate $\la$ and $n\mapsto n-1$
with rate $1_{\{n>0\}}$. Standard results say that this Markov chain is
positive recurrent for $\la<1$, null recurrent for $\la=1$, and transient for
$\la>1$.

For the email process we have just described, in \cite{FS16}, a limit theorem
is proved for the equilibrium distribution of unanswered emails with
priorities just above the critical point $\la_{\rm c}=1$, linking it to the
convex hull of Brownian motion. In \cite{GC09}, a very similar model in
discrete time is studied where priorities are uniformly distributed on $[0,1]$
and in each step, one email is answered and $m\geq 2$ new emails arrive. It is
shown that the critical priority for this model is $p_{\rm c}=1-1/m$ and the
probability that an incoming email has to wait time $t$ before being answered,
given that it is answered at all, decays as $t^{-3/2}$. The proof uses a
mapping to invasion percolation on a regular tree.

The model in \cite{GC09} was inspired by a similar model in \cite{Bar05} where
the number $N$ of emails in the inbox is fixed. In this model, in each step,
with probability $p$ the email with the highest priority is answered, and with
the remaining probability, a random email from the inbox is answered. After
this, a new email with a new random priority arrives. If one sends first $p\to
1$ and then $t\to\infty$, then the probability that an incoming email is
answered after $t$ time steps is asymptotically of order $t^{-1}$. This was
first shown for an inbox containing $N=2$ emails in \cite{Vaz05} and then for
general $N$ in \cite{Ant09}. In the latter paper, it is also argued that
keeping $0<p<1$ fixed and sending first $N\to\infty$ and then $t\to\infty$,
the probability of waiting a time $t$ for an answer decays like $t^{-3/2}$. In
\cite{GC07}, a mapping of the model in \cite{Bar05} to invasion percolation is
described.

\subsubsection*{The Bak-Sneppen model}

The Bak-Sneppen model, introduced in \cite{BS93}, is one of the best-known
models believed to exhibit self-organized criticality. It is a Markov chain
$(X_k)_{k\geq 0}$ with state space $[0,1]^N$. Writing
$X_k=\big(X_k(0),\ldots,X_k(N-1)\big)$, one interprets $X_k(i)\in[0,1]$ as the
fitness of species $i$ at time $k$. One thinks of the species as being
situated on a ring, where $i-1$ and $i+1$ (calculated modulo $N$) are the
neighbors of $i$. Initially, all fitnesses are independent and uniformly
distributed. In each step, the species $i$ with the lowest fitness is
selected, and $X_{i-1},X_i$, and $X_{i+1}$ are all replaced by new,
independent and uniformly chosen fitnesses. It is believed that there exists a
critical value $f_{\rm c}\approx 0.6672(2)$ such that for large $N$, in
equilibrium, the fitness $X(0)$ is approximately uniformly distributed on
$[f_{\rm c},0]$. For $0<f<1$, excursions away from the set $[f,1]^N$ are
called \emph{avalanches}. During an avalanche, to decide what the next move
is, one only needs information about fitnesses below $f$, and in view of this
consecutive avalanches are i.i.d. It is believed that for $f<f_{\rm c}$, the
length of an avalanche and the number of species affected have a limit law as
$N\to\infty$ with exponential tail, but at $f=f_{\rm c}$ one should find
power-law decay signifying self-organized criticality.

The best available rigorous results for the Bak-Sneppen model can be found in
\cite{MZ03,MZ04}. A key tool in the latter paper is the \emph{locking
  thresholds representation}, which shows that at each time, the fitnesses
$X_k(0),\ldots,X_k(N-1)$ can be thought of as being independent and uniformly
distributed in random intervals $[Y_k(0),1],\ldots,[Y_k(N-1),1]$. The paper
defines two critical fitnesses $0<f_{\rm c}\leq f'_{\rm c}<1$, where $f_{\rm
  c}$ resp.\ $f'_{\rm c}$ is the first point where avalanches (in the limit
$N\to\infty$) have infinite expected size, resp.\ are with positive
probability infinite. Assuming that $f_{\rm c}=f'_{\rm c}$, it is shown that
in the limit $N\to\infty$, in equilibrium fitnesses are approximately
uniformly distributed on $[f_{\rm c},1]$.

In \cite{MS12}, a modified Bak-Sneppen model is introduced. Here, if $i$ if
the species with the lowest fitness, then instead of redrawing the fitnesses
of $i-1,i$ and $i+1$, one redraws the fitnesses of $i$ and one random other
species $j$, chosen uniformly from the population. In a sense, this is the
mean-field version of the original Bak-Sneppen model. For this model, it has
been proved in \cite{MS12} that $f_{\rm c}=1/2$. Moreover, the probability
that an avalanche below some chosen fitness level $f$ has a duration longer
than $t$ decays exponentially in $t$ for $f<f_{\rm c}$, but decays as
$t^{-1/2}$ at $f=f_{\rm c}$. This is proved by setting up a coupling with a
branching process, which is subcritical for $f<f_{\rm c}$ and critical for
$f=f_{\rm c}$.

\subsubsection*{The Stigler-Luckock model}

In the Stigler-Luckock model, first introduced in \cite{Sti64} and reinvented
and generalized by Luckock in \cite{Luc03} and again independently reinvented
in \cite{Pla11} and \cite{Yud12a,Yud12b}, traders place buy and sell orders
according to independent Poisson processes. Letting $I=(I_-,I_+)$ denote the
interval of possible prices, buy (resp.\ sell) orders arrive in $A\sub I$ with
intensity $\mu_-(A)$ (resp.\ $\mu_+(A)$), where $\mu_\pm$ are finite measures
on $I$. If a buy order arrives at a price $x$ while the order book already
contains a sell order at a price $y<x$, then the buy order together with the
lowest sell order are immediately removed. Similarly, newly arriving sell
orders at a price $x$ are immediately matched against the best available buy
order at a price $x'>x$, if such an order exists, and stay in the order book
otherwise. Assuming stationarity, Luckock derived a differential equation from
which he was able to calculate the equilibrium distribution of the best buy
and sell offer in the order book. In particular, he was able to calculate two
prices $x_-<x_+$ so that the best sell offer never drops below $x_-$ and the
best buy offer never climbs above $x_+$, with the result that buy orders below
$x_-$ and sell orders above $x_+$ are never matched.

Mathematically, proving existence of the type of stationary distributions that
Luckock postulates remains an open problem. Recent progress was made in
\cite{Swa16}, where for each subinterval $J=(J_-,J_+)\sub I$, a
\emph{restricted model on $J$} is defined in such a way that sell orders that
arrive on the left of $J_-$ can still be matched to the best available buy
order, but if no such order exist, are not written into the order book;
similar rules apply to buy orders that arrive on the right of $J_+$. For such
a restricted model, \cite{Swa16} gives necessary and sufficient conditions for
positive recurrence. In particular, if $x_-,x_+$ are the prices of Luckock,
then it is shown that for any $x_-<J_-<J_+<x_+$, the restricted model is
positive recurrent, while for $J_-<x_-<x_+<J_+$ it is not. Note, however, that
the restricted model on $J$ (as defined above) and the original model
restricted to $J$ evolve in the same way only as long as the best buy and sell
offers remain inside $J$. In view of this, it has so far not been possible to
draw rigorous conclusions about the original model from the behavior of the
restricted model.

The main tool in \cite{Swa16}, which was inspired by Lemma~\ref{L:linfunc}
below which was discovered before \cite{Swa16} was written, are weight
functions that give a weight to buy and sell orders that depends on their
price. By choosing these weight functions in the right way, as solutions to a
suitable differential equation, a Lyapunov function can be constructed. In
\cite{KY16}, a result in the spirit of our Theorem~\ref{T:limsup} is proved
for the Stigler-Luckock model, under certain technical restrictions. An
important tool in their work are comparison lemmas in the spirit of our
Lemmas~\ref{L:monot} and \ref{L:secomp}.

In \cite{PS16}, an extension of the Stigler-Luckock model is introduced with
richer behavior (in particular, a transition between self-organized
criticality and different behavior). See also \cite{Swa16} for references to a
number of related models for the evolution of an order book, which however do
not exhibit self-organized criticality.

In \cite{Pla11,Swa14}, a variation of the Stigler-Luckock model is introduced
where in each step, two traders arrive who want to buy and sell at prices that
are infinitesimally close, but ordered so that they just don't match. This
model can alternatively be interpreted as a model for canyon formation and
motivated the model of the present paper. Numerically, its behavior seems to
be governed by the same prices $x_-,x_+$ as for the Stigler-Luckock model.

\subsubsection*{Other models}

Somewhat similar in spirit to the previous models is also the model
\cite{GMS11}, which is basically a supercitical branching process in which
fitnesses are assigned to the particles, and those killed have the lowest
fitness. We should also mention the branching Brownian motions with strong
selection treated in \cite{Mai16}, which also use a rule of the type ``kill
the lowest particle'' and exhibit self-organized criticality.

We note that in the construction of many of the models discussed so far and in
particular also our model, only the relative order of the points (i.e., their
rank or priority) matters, so replacing the uniform distribution on $[0,1]$ by
any other atomless law on $\R$ yields the same model up to a continuous
transformation of space. Starting from the empty initial state, adding points
one by one, and taking notice only of their relative order, one in effect
constructs after $k$ steps a random permutation of $k$ elements. In view of
this, our quantities of interest may be described as functions of such a
random permutation. This is somewhat reminiscent of the way the authors of
\cite{AD99} use what they call Hammersley's process to study the longest
increasing subsequence of a random permutation. There is an extensive
literature on functions of random permutations, but none of those studied so
far seem relevant for our process.

Also, although our model is an example of a Markov process with
self-reinforcement, it does not seem to have close connections to any of the
classical models with self-reinforcement studied so far, as reviewed in
\cite{Pem07}.

\section{Proofs}

\subsection{Weight functions}

Our main tool for proving Theorems~\ref{T:limsup}--\ref{T:ergod} are linear
functionals that count the number of points in the interval $[0,q]$ weighted
with a function that depends on their position. This method significantly
simplifies methods used in a number of preprints predating the present
paper, which for the interest of the reader are available as \cite{Swa14},
versions 1--3. The method of using weight functions was first discovered for
the present model but has in the mean time already succesfully been applied in
\cite{Swa16} to the Stigler-Luckock model.

As already mentioned in Section~\ref{S:other}, by a simple transformation of
space, we may replace the uniformly distributed random variables $(U_k)_{k\geq
  1}$ by real random variables having any non-atomic distribution. At present,
it will be more convenient to work with exponentially distributed random
variables with mean one, so we transform the unit interval $[0,1]$ into the
closed halfline $[0,\infty]$ with the transformation $q\mapsto
f(q):=-\log(1-q)$ and set $\sig_k:=f(U_k)$ $(k\geq 1)$. Letting
$Y_k:=f(X_k)$ denote the image of $X_k$ under $f$, we then have
\be\label{inductY}
Y_k:=\left\{\ba{ll}
Y_{k-1}\cup\{\sig_k\}\quad&\mbox{if }\sig_k<N_{k-1},\\[5pt]
(Y_{k-1}\cup\{\sig_k\})\beh\{N_{k-1}\}\quad&\mbox{if }\sig_k>N_{k-1}.
\ea\right.\quad(k\geq 1),
\ee
where $N_k:=\min(Y_k\cup\{\infty\})$. We let $\Pc_{\rm fin}[0,\infty]$ denote
the set of all finite subsets of $[0,\infty]$.

\bl[Linear functions of the process]\label{L:linfunc}
For $t\geq 0$, let $\Lin_t:\Pc_{\rm fin}[0,\infty]\to\R$ be defined as
\be\label{Lin}
\Lin_t(Y):=\sum_{s\in Y}e^s1_{[0,t]}(s).
\ee
Then, for the process started in any deterministic initial state $Y_0$, one has
\be\label{Finc}
\E\big[\Lin_t(Y_1)\big]-\Lin_t(Y_0)=t-1_{[0,t]}(N_0).
\ee
\el
\bpro
For any bounded measurable function $w:[0,\infty]\to\R$, we calculate
\bc
\E\big[\Lin_t(Y_1)\big]-F(Y_0)
&=&\dis\int_0^{N_0}w(s)e^{-s}\di s
+\int_{N_0}^\infty\big(w(s)-w(N_0)\big)e^{-s}\di s\\[5pt]
&=&\dis\int_0^\infty w(s)e^{-s}\di s-w(N_0)e^{-N_0}.
\ec
Setting $w(s):=e^s1_{[0,t]}(s)$, we arrive at (\ref{Finc}).
\epro

\subsection{The positive recurrent regime}

For each $t\geq 0$, we let
\be\label{Yt}
Y^{(t)}_k:=Y_k\cap[0,t]\qquad(k\geq 0)
\ee
denote the restriction of the process $Y_k$ to the interval $[0,t]$, which is
itself a Markov chain. Using Lemma~\ref{L:linfunc}, it is easy to show that
for $t<1$, this Markov chain is positive recurrent and ergodic. We let
$\P^\emptyset$ denote the law of the process started in the empty
configuration $Y_0=\emptyset$ and we let
\be\label{tauY}
\tau^{(t)}_\emptyset:=\inf\{k>0:Y^{(t)}_k=\emptyset\}
\ee
denote the first return time of $Y^{(t)}$ to the empty configuration.

\bp[Positive recurrent regime]\label{P:pos}
For $0\leq t<1$, one has
\be\label{return}
\E^\emptyset[\tau^{(t)}_\emptyset]=(1-t)^{-1}.
\ee
Moreover, the process $Y^{(t)}$ has an invariant law $\nu$ on $\Pc_{\rm fin}[0,t]$
such that $\nu(\{\emptyset\})=1-t$ and the process started in an arbitrary
initial law satisfies
\be\label{erg}
\big\|\P[Y^{(t)}_k\in\,\cdot\,]-\nu\big\|\asto{k}0,
\ee
where $\|\,\cdot\,\|$ denotes the total variation norm.
\ep
\bpro
For $t<1$, the function $\Lin_t$ is a Lyapunov function. Using this and the fact
that $\P^\emptyset[Y^{(t)}_1=\emptyset]>0$, which shows that our process is
aperiodic in an appropriate sense, standard results (see, e.g.,
\cite[Prop.~19]{Swa16}) show that $\E^Y[\tau^{(t)}_\emptyset]<\infty$
for any deterministic $Y\in\Pc_{\rm fin}[0,t]$, and there exists
an invariant law $\nu$ such that (\ref{erg}) holds.

Let $(Y^{(t)}_k)_{k\in\Z}$ be the corresponding stationary process. Then
Lemma~\ref{L:linfunc} tells us that
\be
\E[\Lin_t(Y^{(t)}_1)-\Lin_t(Y^{(t)}_0)]=t-\P[Y^{(t)}_0\neq\emptyset].
\ee
By stationarity and Lemma~\ref{L:statinc} below, the left-hand side of this
equation is zero, so the invariant law satisfies $\nu(\{\emptyset\})=1-t$. By
standard renewal arguments,
$\E^\emptyset[\tau^{(t)}_\emptyset]=\nu(\{\emptyset\})^{-1}$, so
(\ref{return}) follows.
\epro

\bl[Stationary increments]\label{L:statinc}
Let $(F(k))_{k\in\Z}$ be a stationary process, and assume that
$\E\big[|F(1)-F(0)|\big]<\infty$. Then $\E\big[F(1)-F(0)\big]=0$.
\el
\bpro
For $M>0$, let $F^M(k):=F(k)$ if $-M\leq F(k)\leq M$ and $F^M(k):=M$ or $-M$
if $F(k)\geq M$ or $F(k)\leq-M$, respectively. By stationarity,
$\E[F^M(1)]=\E[F^M(0)]$ and hence $\E[F^M(1)-F^M(0)]=0$. Letting $M\up\infty$,
using the fact that $|F^M(1)-F^M(0)|\leq|F(1)-F(0)|$ and dominated convergence,
we conclude that $\E[F(1)-F(0)]=0$.
\epro

\subsection{The lower invariant process}

It turns out to be possible to find stationary solutions of the inductive
formula (\ref{inductY}) that are defined for all $k\in\Z$, and this will very
helpful in proving our main theorems. A crucial observation is that solutions
to the inductive formula (\ref{inductY}) are monotone in the starting
configuration.

\bl[First comparison lemma]\label{L:monot}
Let $y$ and $\ti y$ be finite subsets of $[0,1]$ and let $(Y_k)_{k\geq
  0}$ and $(\ti Y_k)_{k\geq 0}$ be defined by the inductive relation
(\ref{inductY}) with $Y_0=y$ and $\ti Y_0=\ti y$. Then $y\sub\ti
y$ implies that $Y_k\sub\ti Y_k$ for all $k\geq 0$.
\el
\bpro
It suffices to show that $Y_{k-1}\sub\ti Y_{k-1}$ implies
$Y_k\sub\ti Y_k$. Adding the point $\sig_k$ to both
$Y_{k-1}$ and $\ti Y_{k-1}$ obviously preserves the order
of inclusion, as does simultaneously removing the minimal elements
$N_{k-1}$ from $Y_{k-1}$ and $\ti N_{k-1}$ from $\ti Y_{k-1}$. Since
$Y_k\sub\ti Y_k$ we have $N_{k-1}\geq\ti N_{k-1}$ and it may happen that
$\ti N_{k-1}<\sig_k\leq N_{k-1}$, in which case we remove
$\ti N_{k-1}$ from $\ti Y_{k-1}$ but not $N_{k-1}$ from
$Y_{k-1}$, but in this case $\ti N_{k-1}$ is not an element
of $Y_{k-1}$ so again the order is preserved.
\epro

We will be interested in stationary solutions to the inductive relation
(\ref{inductY}). To this aim, we consider a two-way infinite sequence
$(\sig_k)_{k\in\Z}$ of i.i.d.\ exponentially distributed random variables with
mean one. For each $m\in\Z$, we let $(Y_{m,k})_{k\geq m}$ denote the solution
to the inductive relation (\ref{inductY}) started in $Y_{m,m}:=\emptyset$.
Since $Y_{m-1,m}\supset\emptyset=Y_{m,m}$, we see by Lemma~\ref{L:monot}, that
$Y_{m-1,k}\supset Y_{m,k}$ for all $k\geq m$, so there exists a collection
$(\Yo_k)_{k\in\Z}$ of countable subsets of $\half$ such that
\be\label{lowinv}
Y_{m,k}\up \Yo_k\quad\mbox{as }m\down-\infty.
\ee
We call the limit process $(\Yo_k)_{k\in\Z}$ from (\ref{lowinv}) the \emph{lower
  invariant process}. The following proposition is the main result of the
present subsection.

\bp[Lower invariant process]\label{P:lowinv}
For all $k\in\Z$, one has
\be
\big|\Yo_k\cap[0,t]\big|\left\{\ba{ll}
<\infty\quad{\rm a.s.}\quad&\mbox{if }t\in[0,1),\\[5pt]
=\infty\quad{\rm a.s.}\quad&\mbox{if }t\in[1,\infty).\ea\right.
\ee
The set $\Yo_k$ a.s.\ has a minimal element
$\No_k:=\min(\Yo_k)$, whose distribution is given by
\be\label{Nt}
\P[\No_k<t]=t\wedge 1\qquad(k\in\Z,\ t\geq 0).
\ee
Moreover, the process $(\Yo_k)_{k\in\Z}$ solves the inductive relation
(\ref{inductY}) for all $k\in\Z$.
\ep
\bpro
Proposition~\ref{P:pos} tells us that for each $t\in[0,1)$ and fixed $k\in\Z$,
the random variables $Y_{m,k}$ converge in distribution as $m\down-\infty$
to a limit with law $\nu$ satisfying $\nu(\{\emptyset\})=1-t$. This shows
that $\Yo_k\cap[0,t]$ is a.s.\ finite for each $t<1$ and
$\P\big[\Yo_k\cap[0,t]=\emptyset\big]=1-t$. It follows that $\Yo_k\cap[0,1)$ is a
locally finite subset of $[0,1)$ and
\be
\P\big[\Yo_k\cap[0,1)=\emptyset\big]
=\lim_{t\up 1}\P\big[\Yo_k\cap[0,t]=\emptyset\big]=0.
\ee
If $\P\big[\big|\Yo_k\cap[0,1)\big|\leq N\big]$ were positive for some
$N<\infty$, then it is easy to see that $\P[\Yo_{k+N}=\emptyset]$ would also
be positive, so by stationarity we conclude that $\Yo_k\cap[0,1)$ is
a.s.\ infinite. In particular, this implies that $\Yo_k\cap[0,t]=\infty$
a.s.\ for each $t\geq 1$. Moreover, (\ref{Nt}) now follows from the fact that
$\P\big[\Yo_k\cap[0,t]=\emptyset\big]=1-t$.

Our arguments so far show that almost surely, $\No_{k-1}<1$ and
$\Yo_{k-1}\cap[0,t]$ is a finite set for all $t<1$. Choose $\No_{k-1}<t<1$.
Then a.s.\ there exists an $M>-\infty$ such that
$Y_{m,k-1}\cap[0,t]=\Yo_{k-1}\cap[0,t]$ for all $m\leq M$. Using this and the
fact that $Y_{m,k-1}$ and $Y_{m,k}$ are related by the inductive relation
(\ref{inductY}), we see that also $\Yo_{k-1}$ and $\Yo_k$ are related as in
(\ref{inductY}).
\epro

\subsection{Proof of the theorems}

Using the transformation $\sig_k=-\log(1-U_k)$, Theorem~\ref{T:reccur} can
equivalently be formulated for the process $Y$ from (\ref{inductY}) as
follows.

\bt[Recurrence versus transience]\label{T:Yreccur}
Let $Y^{(t)}$ be the restricted process from (\ref{Yt}) and let
$\tau^{(t)}_\emptyset$ as in (\ref{tauY}) denote the first return time of
$Y^{(t)}$ to the empty configuration. Then
\be\ba{ll}\label{Yretn}
\dis\E^\emptyset[\tau^{(t)}_\emptyset]=(1-t)^{-1}\qquad&\dis(t<1),\\[5pt]
\dis\E^\emptyset[\tau^{(t)}_\emptyset]=\infty\quand
\P^\emptyset[\tau^{(t)}_\emptyset<\infty]=1\qquad&\dis(t=1),\\[5pt]
\dis\P^\emptyset[\tau^{(t)}_\emptyset=\infty]>0\qquad&\dis(t>1).
\ec
\et
\bpro
The case $t<1$ has already been proved in Proposition~\ref{P:pos}.
Clearly $t\leq t'$ implies $\tau^{(t)}_\emptyset\leq\tau^{(t')}_\emptyset$
a.s., so the formula for $t<1$ implies that
$\E^\emptyset[\tau^{(t)}_\emptyset]=\infty$ for all $t\geq 1$.

For $t=1$, Lemma~\ref{L:linfunc} implies that
\be
M_k:=\Lin_1(Y_{k\wedge\tau^{(1)}_\emptyset})\qquad(k\geq 0)
\ee
is a nonnegative martingale, so the a.s.\ limit
$M_\infty:=\lim_{k\to\infty}M_k$ exists. It is easy to see that $M_k$ cannot
converge to a positive limit, so we conclude that $M_\infty=0$ a.s.\ and hence
$\tau^{(1)}_\emptyset<\infty$ a.s.

For $t>1$, we observe that for any $0\leq s<t$,
\be
\big|Y_n\cap(s,t]\big|\geq\sum_{k=1}^n1_{\txt\{s<\sig_k<t\}}
-\sum_{k=1}^n1_{\txt\{Y_{k-1}\cap[0,s]=\emptyset\}}.
\ee
By the strong law of large numbers and Proposition~\ref{P:pos}
\be
n^{-1}\sum_{k=1}^n1_{\txt\{s<\sig_k<t\}}\asto{n}(e^{-s}-e^{-t})
\quand
n^{-1}\sum_{k=1}^n1_{\txt\{Y_{k-1}\cap[0,s]=\emptyset\}}\asto{n}1-s\quad{\rm a.s.}
\ee
Choosing $s$ close enough to 1 such that $1-s<e^{-s}-e^{-t}$, we see that
$\big|Y_n\cap(s,t]\big|\to\infty$ a.s., and hence 
$\P^\emptyset[\tau^{(t)}_\emptyset=\infty]>0$.
\epro

\bpro[of Theorem~\ref{T:limsup}]
It follows from Proposition~\ref{P:pos} that regardless of the initial state,
$\limsup_{k\to\infty}N_k\geq 1$. On the other hand, in the proof of
Theorem~\ref{T:Yreccur} we have seen that $|Y_n\cap[0,t]|\to\infty$ a.s.\ for
any $t>1$, proving that $\limsup_{k\to\infty}N_k\leq 1$. Translating these
results to the process $X$ yields Theorem~\ref{T:limsup}.
\epro

\bpro[of Theorem~\ref{T:ergod}]
This follows from results that have already been proved for the process $Y$
from (\ref{inductY}). The ergodic statement (\ref{ergod}) follows from
Proposition~\ref{P:pos}, and the fact that $X_\infty$ is an infinite, but
locally finite subset of $[0,p_{\rm c})$ follows from Proposition~\ref{P:lowinv}.
\epro

\subsection{Some concluding remarks}\label{S:conclude}

For the process that is the subject of the present paper, two open problems
seem worth investigating. For the transformed process $Y$ from
(\ref{inductY}), one would like to know if it is true, as numerical
simulations suggest, that at the critical point $t_{\rm c}:=1$, the tail of
the return probability decays as $\P^\emptyset[\tau^{(t_{\rm c})}\geq k]\sim
k^{-1/2}$. Second, one would like to know the asymptotic shape of the locally
finite point set $\Yo_k\cap[0,1)$ of Proposition~\ref{P:lowinv} near the
  critical point $t_{\rm c}=1$.

Both these problems have been resolved for the email model described in
Section~\ref{S:other}. In particular, in \cite{FS16}, for the email model, it
has been shown that the (random) distribution function of the lower invariant
process, properly rescaled, converges to a limit as one approaches the critical
point, and the limit can be expressed in terms of the convex hull of Brownian
motion. The proof of this fact essentially uses that for the email process,
the number of emails in the inbox above a certain priority is a random walk
with reflection at the origin, which properly rescaled converges to a
(drifted) Brownian motion.

To prove analogue results for the present, more complicated canyon model, it
seems essential to understand the point-counting process
\be\label{Fdef}
F_t(k):=\big|Y_k\cap[0,t]\big|\qquad\big(k\in\N,\ t\in\half\big),
\ee
and similarly $\Fo_t(k)$, which is defined in terms of the lower invariant
process $\Yo$. Although these are not Markov processes, one possible guess is
that for $t<1$ close to one, suitable rescaled, they converge to drifted
Brownian motions with reflection at the origin. Alternatively, defining
$\Lin_t$ as in (\ref{Lin}), one may also look at the processes $\Lin_t(Y_k)$
and $\Lin_t(\Yo_k)$, which count particles in $[0,t]$ weighted with the
function $e^x$. In view of Lemma~\ref{L:linfunc}, one has good control over
the compensator of these processes and hence convergence to reflected Brownian
motion may be easier to prove.

We mention two technical facts that were used in preprints preceding the
present paper (\cite{Swa14}, versions 1--3) and that may still be of interest.

First, due to the memoryless property of the exponential distribution of the
random variables $\sig_k$, it is not hard to check that the function valued
process $(F_t)_{t\geq 0}$, with $F_t=(F_t(k))_{k\geq 0}$ as in (\ref{Fdef}),
is a continuous-time Markov process, where the parameter $t$ plays the role of
time. Indeed, at the time $t=\sig_k$, the function $F_t$ changes as
\be\label{Fjump}
F_t(k')=\left\{\ba{ll}
F_{t-}(k')+1\quad&\mbox{if }k\leq k'<\kappa_t(k),\\[5pt]
F_{t-}(k')\quad&\mbox{otherwise}\ea\right.
\qquad(k'\geq 0),
\ee
where $F_{t-}$ denotes the state immediately prior to time $t$ and
\be\label{kappa}
\kappa_t(k):=\inf\{k'>k:F_{t-}(k'-1)=0=F_{t-}(k')\},
\ee
with the convention that $\inf\emptyset:=\infty$. Similarly, $(\Fo_t)_{0\leq
  t<1}$, with $\Fo_t=(\Fo_t(k))_{k\in\Z}$ is also a Markov process, which,
however, is only well-defined until time one. Also $\Lin_t(Y_k)$
and $\Lin_t(\Yo_k)$ evolve as a function of $t$ in a Markovian way, 
but in this case, in (\ref{Fjump}), the term $+1$ needs to be replaced by
$+e^t$.

A second technical fact that is perhaps of interest is a second comparison
lemma, similar to Lemma~\ref{L:monot}, but using a different order. Roughly
speaking, it says that for $Y_k\cap[0,t]$ to avoid becoming the empty set, it
is good to have many particles that are situated as far as possible to right
in the interval $[0,t]$.

\bl[Second comparison lemma]\label{L:secomp}
For each $0\leq s\leq t$ and finite $y\sub\half$, let
$F^y_{s,t}(k):=\big|Y^y_k\cap[s,t]\big|$ $(k\geq 0)$. Fix $t>0$ and let
$x,y\sub\half$ be finite. Then
\be\label{secomp}
F^x_{s,t}(0)\leq F^y_{s,t}(0)\ \forall s\in[0,t]
\quad\mbox{implies}\quad
F^x_{s,t}(k)\leq F^y_{s,t}(k)\ \forall s\in[0,t],\ k\geq 0.
\ee
\el
\bpro
It suffices to prove (\ref{secomp}) for $k=1$; the general statement follows
by induction. Without loss of generality, we may also assume that $x$ and $y$
are subsets of $[0,t]$. Order the elements of $x$ and $y$ as
$x=\{x_1,\ldots,x_n\}$ and $y=\{y_1,\ldots,y_m\}$ with $x_n<\cdots<x_1$ (in
this order!) and $y_m<\cdots<y_1$.  Then the assumption that $F^x_{s,t}(0)\leq
F^y_{s,t}(0)\ \forall s\in[0,t]$ is equivalent to the statement that $m\geq n$
and $x_i\leq y_i$ for all $i=1,\ldots,n$. We must show that we can order the
elements of $\ti x:=Y^x_1\cap[0,t]$ and $\ti y:=Y^y_1\cap[0,t]$ in the same
way. We distinguish three different cases.

Case~I: $\sig_1<x_n$. In this case, no points are removed from $x$ while $\ti
x_{n+1}:=\sig_1$ is added as the $(n+1)$-th element. Since $x_n\leq y_n$, the
elements $y_1,\ldots,y_n$ remain unchanged while $\ti y_{n+1}$ is the maximal
element of $\{\sig_1\}\cup\{y_m,\ldots,y_{n+1}\}$, which lies on the right of
$\ti x_{n+1}=\sig_1$.

Case~II: $x_n<\sig_1<t$. In this case, $x_n$ is removed from $x$ and there
exist $1\leq n'\leq n$ and $n'\leq m'\leq m+1$ such that $\sig_1$ is inserted
into $x$ between the $n'$-th and $(n'-1)$-th element and into $y$ between the
$m'$-th and $(m'-1)$-th element, where we allow for the cases that $n'=1$
($\sig_1$ is added at the right end of $x$ and possibly also of $y$) and
$m'=m+1$ ($\sig_1$ is added at the left end of $y$). The elements of the new
sets $\ti x$ and $\ti y$, ordered from low to high, are now
\be\ba{r@{\,}r@{\,}r@{\,}c@{\,}l}
\{x_{n-1},\ldots,x_{m'},&x_{m'-1},\ldots,x_{n'},\sig_1,&x_{n'-1},\ldots,x_1\}
&=&\ti x,\\[5pt]
\{y_{m-1},\ldots,y_n,
y_{n-1},\ldots,y_{m'},&\sig_1,y_{m'-1},\ldots,y_{n'},&y_{n'-1},\ldots,y_1\}
&=&\ti y.
\ec
Here $x_{n-1},\ldots,x_{m'}$ lie on the left of $y_{n-1},\ldots,y_{m'}$,
and likewise $x_{n'-1},\ldots,x_1$ lie on the left of $y_{n'-1},\ldots,y_1$,
respectively. Since moreover
\be
x_{m'-1}<\cdots<x_{n'}<\sig_1<y_{m'-1}<\cdots<y_{n'},
\ee
these elements are ordered in the right way too.

Case~III: $t<\sig_1$. In this case, the lowest elements of $x$ and $y$ are
removed while no new elements are added, which obviously also preserves the
order.
\epro

\subsection*{Acknowledgements}

The main method used in the proofs was discovered in a joint discussion with
Marco Formentin. I would further like to thank Persi Diaconis, Patrik Ferrari,
Ronald Meester, Martin Ondrej\'at, Frank Redig, Jan Seidler, Martin
\v{S}m\'id, and Rongfeng Sun for useful discussions and in particular Jan van
Neerven, Marco Formentin, and Martin \v{S}m\'id for drawing my attention to
the Bak-Sneppen model, Barab\'asi's queueing system, and the Stigler-Luckock
model, respectively.


\begin{thebibliography}{DGW05}

\bibitem[AD99]{AD99}
D.~Aldous and P.~Diaconis.
Longest increasing subsequences: from patience sorting to the
Baik-Deift-Johansson theorem.
\emph{Bull.\ Am.\ Math.\ Soc., New Ser.}~36(4) (1999), 413--432.

\bibitem[Ant09]{Ant09}
C.~Anteneodo.
Exact results for the Barab\'asi queuing model.
\emph{Phys.\ Rev.\ E}~80 (2009), 41131--41137.

\bibitem[Bar05]{Bar05}
A.-L.~Barab\'asi.
The origin of bursts and heavy tails in human dynamics.
\emph{Nature}~435 (2005), 207--211.

\bibitem[BS93]{BS93}
P.~Bak and K.~Sneppen.
Punctuated equilibrium and criticality in a simple model of evolution.
\emph{Phys.\ Rev.\ Lett.}~74 (1993), 4083--4086.

\bibitem[FS16]{FS16}
M.~Formentin and J.M.~Swart.
The limiting shape of a full mailbox.
\emph{ALEA}~13(2) (2016), 1151--1164.

\bibitem[GC07]{GC07}
A.~Gabrielli and G.~Caldarelli.
Invasion Percolation and Critical Transient in the Barab\'asi Model of Human
Dynamics.
\emph{Phys.\ Rev.\ Lett.}~98 (2007), 208701--208704.

\bibitem[GC09]{GC09}
A.~Gabrielli and G.~Caldarelli.
Invasion percolation and the time scaling behavior of a queuing model of human
dynamics.
\emph{J.\ Stat.\ Mech.} (2009), P02046 (10 pages).

\bibitem[GMS11]{GMS11}
H.~Guiol, F.P.~Machado, and R.B.~Schinazi.
A stochastic model of evolution.
\emph{Markov Process.\ Relat.\ Fields}~17(2) (2011), 253--258.

\bibitem[KY16]{KY16}
F.~Kelly and E.~Yudovina.
A Markov model of a limit order book: thresholds, recurrence,
and trading strategies.
Preprint (2016), 24 pages, ArXiv:1504.00579v3.

\bibitem[Luc03]{Luc03}
H.~Luckock.
A steady-state model of the continuous double auction.
\emph{Quantitative Finance}~3(5) (2003), 385--404.

\bibitem[Mai16]{Mai16}
P.~Maillard.
Speed and fluctuations of N-particle branching Brownian motion with spatial
selection.
\emph{Probab.\ Theory Relat.\ Fields}~166(3) (2016), 1061--1173.

\bibitem[MS12]{MS12}
R.~Meester and A.~Sarkar.
Rigorous self-organised criticality in the modified Bak-Sneppen model.
\emph{J.\ Stat.\ Phys.}~149 (2012), 964--968.

\bibitem[MZ03]{MZ03}
R.~Meester and D.~Znamenski.
Limit behavior of the Bak-Sneppen evolution model.
\emph{Ann.\ Probab.}~31(4) (2003), 1986--2002.

\bibitem[MZ04]{MZ04}
R.~Meester and D.~Znamenski.
Critical thresholds and the limit distribution in the Bak-Sneppen model.
\emph{Commun.\ Math.\ Phys.}~246(1) (2004), 63--86.

\bibitem[Pem07]{Pem07}
R.~Pemantle.
A survey of random processes with reinforcement.
\emph{Probab.\ Surveys}~4 (2007), 1--79.

\bibitem[Pla11]{Pla11}
Jana Pla\v{c}kov\'a.
Shluky volatility a dynamika popt\'avky a nab\'idky.
Master Thesis, MFF, Charles University Prague, 2011.

\bibitem[PS16]{PS16}
V.~Per\v{z}ina and J.M.~Swart.
How many market makers does a market need?
Preprint (2016), 14 pages, ArXiv:1612.00981.

\bibitem[Sti64]{Sti64}
G.J.~Stigler.
Public Regulation of the Securities Markets
\emph{The Journal of Business}~37(2) (1964), 117--142.

\bibitem[Swa14]{Swa14}
J.M.~Swart.
A self-reinforced model for canyon formation.
Preprint (2014) ArXiv:1405.3609v1.

\bibitem[Swa16]{Swa16}
J.M.~Swart.
Rigorous results for the Stigler-Luckock model for the evolution of an order
book.
Preprint (2016) ArXiv:1605.01551.

\bibitem[Tok15]{Tok15}
I.M.~Toke.
The order book as a queueing system: average depth and influence of the size
of limit orders.
\emph{Quantitative Finance}~15(5) (2015), 795--808.

\bibitem[Vaz05]{Vaz05}
A.~V\'azquez.
Exact results for the Barab\'asi model of human dynamics.
\emph{Phys. Rev. Lett.}~95 (2005), 248701--248704.

\bibitem[Yud12a]{Yud12a}
E.~Yudovina.
A simple model of a limit order book.
Preprint (2012), 27 pages, ArXiv:1205.7017v2

\bibitem[Yud12b]{Yud12b}
E.~Yudovina.
Collaborating Queues: Large Service Network and a Limit Order Book.
Ph.D.\ thesis, University of Cambridge, 2012.

\end{thebibliography}
\end{document}